\documentclass[11pt]{article} 
\usepackage{amsmath}
\usepackage{makeidx}
\usepackage{graphicx}
\usepackage{epsfig}
\DeclareMathOperator*{\osc}{Osc}
\newtheorem{theorem}{Theorem}

\newtheorem{definition}{Definition}
\newtheorem{lemma}{Lemma}

\def\intl{\int\kern-9pt\hbox{$\backslash$}}
\begin{document}
\textwidth 5.5 truein
\title{\vspace{-3cm} Geometric De Giorgi Theory}
\author{Lihe Wang}
\date{}
\maketitle
The one of the most fundamental and influential theories for elliptic equations  shall be the
De Giorgi-Nash-Moser theorem on the continuity of weak solutions to uniformly elliptic equations.
The theorem was first discovered by E. De Giorgi \cite{DG} in 1957  and independently by J. Nash \cite{Nash} in 1958.    
Later J. Moser \cite{Moser} in 1960 gave a proof using  the powers of the solution and John-Nirenberg inequality.  Late on, Morrey \cite{Morrey} used logarithm of the solution to simplify one step of the argument.

De Giorgi used the theorem to solve Hilbert's 19th problem and his  proof has far  reaching influence  on the  whole field of the analysis of partial differential equations.  Since then these proofs have fascinated many mathematicians in their ingenuity as well as the delicate iterations.

The key of our approach is that the theorem follows from two new simple lemmas \ref{lambda} and \ref{key2} and an iteration of these implies the theorem. The novelty is that each lemma have clear geometrical interpretation and the computation is also geometrical and intuitive. 

We use the standard notions. $x=(x_1, x_2, \dots, x_n)$ denote a point in the $n$-dimensional Euclidean space.
$B_r$ is the ball with center $0$. $D_i u=\frac {\partial u }{\partial x_i}$ and $Du=(D_1u, D_2u, \dots, D_n u)$. If $E$ is a measurable set, $|E|$ is its Lebesgue measure
while if $E$ is of $(n-1)$ dimensional and if no confusion arises, we also use $|E|$ for its the $(n-1)$-dimensional measure. We also use the convention that repeated indexes are summed from 1 to $n$.

The plan of this paper is that we will prove the maximum bound  in Section 1 and oscillation estimate in section 2.

\vglue 1pc
\section{Geometric Local $L^\infty$  Theory} 
\smallskip
We will consider equations in divergent form as
\begin{equation}
Lu=-D_j(a_{ij}(x)D_iu)=0\label{Lu}
\end{equation}
in this paper.  We assume the coefficients satisfy the {\it Uniformly Ellipticity Condition}  if there are two constants
$0<\lambda\leq \Lambda<\infty,$ which are called ellipticity constants, so that 
 \begin{align}
\lambda|\xi|^2\le &a_{ij}(x)\xi_i\xi_j\le \Lambda|\xi|^2, \mbox{ for all } \xi\in R^n,\label{elliptic}\\
|a_{ij}|\le & \Lambda.
\end{align}

\begin{definition} 
\index{weak solution!subsolution} \index{weak
solution!supersolution} We say $u\in H^1(B_r)$ is a
subsolution (solution) of (\ref{Lu}) in $B_r$ write as
$$
Lu=-D_j(a_{ij}(x)D_iu)\leq(=) 0
$$
if for any nonnegative $ v\in H_0^1(B_r)$, we have
$$
\int a_{ij}D_iuD_jvdx\leq ( =) 0.
$$
\end{definition} 

The following two lemmas are standard and we include them for the sake of completeness.

\begin{lemma} \label {1} 
If $u$ is a subsolution function of (\ref{Lu}),
so is $u^+$.
\end{lemma} 
\vglue 1pc
{\bf Proof.} Let $G(s)=s^+$ and for $\epsilon>0$, let $G_\epsilon$ be the convex function as
$$
G_\epsilon(s)=\left\{\begin{array}{ll}
	0, &\mbox{ if } s\leq 0;\\
	\frac 1{2\epsilon}s^2, &\mbox{ if } 0<s\leq \epsilon;\\
	s+\frac{\epsilon}2, &\mbox{else}.
	\end{array}\right.
	$$
	
We observe,  that $G_\epsilon\in C^{1,1}(R^1)$ and increasing and convex as $G_\epsilon^\prime (s)\ge 0$,
$G_\epsilon^{\prime\prime}(s)\ge 0$.  Thus
we have that $G_\epsilon(u)\in H^1$ and $D_i(G_\epsilon(u))=G_\epsilon^\prime(u)D_i u$ from
Theorem 2.2.2, page 52 in \cite{Zie}. Now, for any nonegative $v\in H_0^1 (B_1)$,
 we have
$$
\begin{array}{rcl} 
\int a_{ij}(G_\epsilon(u))_iv_jdx&=&\int a_{ij}G_\epsilon^\prime(u)D_iuD_jvdx \\
&=&\int a_{ij}D_iuD_j(G_\epsilon^\prime(u)v)dx-\int a_{ij}G_\epsilon^{\prime\prime}
(u)D_iuD_jv dx \\
&\le& 0
\end{array}
$$
\noindent by using $u$ is a subsolution.\smallskip

\noindent The lemma follows by letting $\epsilon\to 0$ using the dominated convergence theorem.
\vglue 1pc

We need the following lemma, which is called local energy estimate
or Caccioppoli inequality. \index{Caccioppoli inequality!see {local
energy estimate}}

\begin{lemma} \label{energy} 
 {\rm{(Energy Estimate for Subsolutions)}}
There are  constants $C_1$ and $C_2$ depending on the ellipticity constants and the dimension that for any nonnegative subsolution $u$ of (\ref{Lu}) in $B_1$ the following statements are true:  

(a) For any Lipschitz function $\eta$
with zero boundary value on  $\partial B_1$, we have

\begin{equation}
\int\eta^2|D u|^2\le C_1\int |D\eta|^2u^2\ dx\ .\label{energy2}
\end{equation}

(b) If $u$ is a nonnegative subsolution in $B_1$ of (\ref{Lu}),
then 
\begin{equation}\label{energy1}
\int_{B_\frac 12}|D u|^2\le 4C_1\int_{B_1} u^2.
\end{equation}

(c) 
If $u$ is a nonnegative subsolution in $B_r$ of (\ref{Lu}),
then 
\begin{equation}\label{energyr}
r^2\int_{B_\frac r2}|D u|^2\le 4C_1\int_{B_r} u^2.
\end{equation}

(d) If $u$ is a nonnegative subsolution in $B_1$ of (\ref{Lu}), then
\begin{equation}\label{gain}
\left(\int_{B_\frac 12} u^{\frac{2n}{n-2}}\right)^{\frac {n-2}{2n}}\leq C_2 \left(\int_{B_1} u^2\right)^\frac 12
\end{equation}
( if $n=2$,  we replace $\frac {2n}{n-2}$ by a number bigger than 2).
\end{lemma} 

\vglue 1pc
{\bf Proof:} 
(a) Taking $\eta^2u$ as a test function,
$$
0\ge \int D_j(\eta^2u) a_{ij}D_i u=
\int \eta^2D_jua_{ij}D_iu+2
\int u\eta D_j\eta a_{ij}D_iu\ .
$$

\noindent By the ellipticity condition,
\begin{align*}
\lambda\int\eta^2|D u|^2&\le\int\eta^2a_{ij}D_iuD_ju\\
&\le -2\int\eta D_i\eta a_{ij}uD_ju\\
&\le  2n\Lambda\int\eta|D\eta||D u|u\\
&\le \frac 12 \lambda\int\eta^2|D u|^2+2\frac{n^2\Lambda^2}
{\lambda}\int|D\eta|^2u^2.
\end{align*}

\noindent (\ref{energy2}) follows immediately by canceling  $\frac 12\lambda
\int\eta^2|Du|^2$ for $C_1=\frac {4n^2\Lambda^2}{\lambda^2}\geq 1$. 

(b) 
Now we take 
$$
\eta(x)=
\left\{\begin{array}{ll} 1,& \mbox{ if } x\in B_\frac 12,\\
2({1-|x|}),& \mbox{ if } x\in B_1\backslash B_\frac 12.\end{array}\right.
$$

Therefore $|D \eta|\leq 2$ and thus
$$
\int_{B_\frac 12}|D u|^2\le 4C_1\int_{B_1} u^2.
$$

(c) This follows the scaling, that is applying   (b) on $u(rx)$.

(d) From (a), we have
\begin{equation}
\int|D(\eta  u)|^2\le 2 C_1\int |D\eta|^2u^2\ dx\ .\label{energy3}
\end{equation}

Let  $S$ be the embedding constant of $H^1_0(B_1)\subset L^\frac {2n}{n-2}(B_1)$,
as for any $v\in H^1_0(B_1)$,
$$\left(\int_{B_1}\;\; v^{\frac{2n}{n-2}}\right)^{\frac {n-2}{2n}}\leq S \left(\int_{B_1}|D v|^2\right)^\frac 12,
$$
and apply this to $\eta u\in H^1_0(B_1)$ with the $\eta$ as in the proof of (b) above,

\begin{align*}
\left(\int_{B_\frac 12}u^{\frac{2n}{n-2}}\right)^{\frac {n-2}{2n}}\leq &
\left(\int_{B_1} |\eta u|^{\frac{2n}{n-2}}\right)^{\frac {n-2}{2n}}\\
\leq &S \left(\int_{B_1}|D\eta u|^2\right)^\frac 12\\
\leq& 2S\sqrt{2C_1}\left(\int_{B_1}u^2\right)^\frac 12.
\end{align*}
This proves that (d) is true for $C_2=2^{\frac 32}SC_1^\frac 12$.

\medskip

We prove the first iteration lemma which yields the local maximum bound.

\vglue 1pc
\begin{lemma}\label{lambda}
(a) There is an $A>1$ depending only on the dimension and ellipticity constants so that if  $u$ is a nonnegative subsolution of (\ref{Lu}) in $B_1$ with its  $L^2$ average  over  $B_1$  defined as
$$
m=\left(\int_{B_1}\kern-16pt\hbox{$\backslash$} \;\; u^2 dx\right)^\frac 12,
$$
then the  $L^2$ average over $B_\frac 12$ of $v=(u(x)-Am)^+$  has estimate:
\begin{equation}
\left(\int_{B_\frac 12}\kern-18pt\hbox{$\backslash$}\;\; v^2 dx\right)^\frac 12\leq \frac 12\left(\int_{B_1}\kern-15pt\hbox{$\backslash$} \;\; u^2 dx\right)^\frac 12.
\end{equation}

(b)  If  $u$ is a nonnegative subsolution of (\ref{Lu}) in $B_r$ with its  $L^2$ average $m$ over  $B_r$  as
$$
m=\left(\int_{B_r}\kern-16pt\hbox{$\backslash$} \;\; u^2 dx\right)^\frac 12,
$$
then for $v=(u(x)-Am)^+$, we have 
 \begin{equation}\label{1key}
\left(\int_{B_\frac r2}\kern-18pt\hbox{$\backslash$}\;\; v^2 dx\right)^\frac 12\leq \frac 12\left(\int_{B_r}\kern-15pt\hbox{$\backslash$} \;\; u^2 dx\right)^\frac 12.
\end{equation}
\end{lemma}

\medskip

{\bf Remark 1.} The lemma and its proof are using the geometric interpretation of the estimate (\ref{gain}) from $L^2$ to $L^{\frac {2n}{n-2}}$. The essence  of this estimate is a faster decay of the distribution  of the levelsets:
$$
\{x: u(x)\geq t\}.
$$
 However,  the constant $C_2$ makes it not significant at lower levelsets so we have to go higher enough level, i.e, as $Am$  in the lemma, to exploit the full power of this decay.

\medskip

\noindent{\bf Proof.} We prove (a) only since (b) follows from (a) by scaling. First we observe $v\leq u$ and from  (\ref{gain}) we have that
$$
\left(\int_{B_\frac 12} v^{\frac{2n}{n-2}}\right)^{\frac {n-2}{2n}}\leq \left(\int_{B_\frac12} u^{\frac{2n}{n-2}}\right)^{\frac {n-2}{2n}}\\
\leq C_2 \left(\int_{B_1} u^2\right)^\frac 12=C_2|B_1|^\frac 12m,
$$
and from Chebyshev inequality,
\begin{align*}
|\{x\in B_\frac 12: v\geq 0\}|&=|\{x\in B_\frac 12: u\geq Am\}|\\
&\leq \frac 1{(Am)^2}\int_{B_\frac 12} u^2 dx\\
&\leq \frac 1{(Am)^2}\int_{B_1} u^2 dx\\
&=\frac {|B_1|}{A^2}.
\end{align*}

Therefore, we use the H\"older inequality combining these two inequalities to obtain,
\begin{align}\nonumber
\int_{B_\frac 12}v^2 dx &\leq \left|\{x\in B_\frac 12: v\geq 0\}\right|^\frac 2n\left(\displaystyle\int_{B_\frac 12} v^{\frac{2n}{n-2}}\right)^{\frac {n-2}{n}}\\[.2in]
&\nonumber\leq \frac {|B_1|^\frac 2n}{A^{\frac 4 n}}\cdot C_2^2|B_1|m^2\\[.2in]
&= \frac 1{4}|B_\frac 12|m^2,\label{iden}
\end{align}
if we take 
\begin{equation}A=\frac{2^{\frac{n(n+7)}4}n^\frac {n}{2}S^{\frac n2}}{|B_1|^{\frac {n-2}4}}\left(\frac \Lambda\lambda\right)^\frac n2.\end{equation}
  This completes the proof.
\vglue 1pc

\begin{lemma}\label {point-infty}
Suppose $u$ is nonnegative subsolution of (\ref{Lu}) in $B_1$ and suppose $0$ is a Lebesgue point for $u$.
Then  with the constant $A$ defined in Lemma \ref{lambda}, we have
$$u(0)\leq 2A\left(\int_{B_1}\kern-16pt\hbox{$\backslash$} \;\; u^2 dx\right)^\frac 12.$$
\end{lemma}

{\bf Proof.} Let $v_0=u$ and  for $k=1,2, \dots,$ let $$
\begin{array}{rl}
v_k &=(v_{k-1}-\frac {Am_{k-1}}{2^k})^+,\\
 m_{k-1}&=\left(\displaystyle\int_{B_{\frac 1{2^{k-1}}}}\kern-32pt\hbox{$\backslash$}\;\; \;\; v_{k-1}^2 dx\right)^\frac 12.
\end{array}
$$

Now we prove the following by induction on $k=0,1,\dots$,
\begin{equation}\label{key1}
\left(\int_{B_{\frac 1{2^{k+1}}}}\kern-32pt\hbox{$\backslash$}\;\;\;\; v_{k+1}^2 dx\right)^\frac 12\leq \frac 12\left(\int_{B_{\frac 1{2^{k}}}}\kern-22pt\hbox{$\backslash$}\;\; \;\; v_k^2 dx\right)^\frac 12.
\end{equation}

First we notice that for $k=0$ is exactly previous lemma. Now assume $k$. 
Now  for the $k+1$ step,  it follows apply Lemma (\ref{1key}) for  on $v_k$ in $B_\frac 1{2^k}$ to obtain the estimate
$$ 
\left(\int_{B_{\frac 1{2^{k+1}}}}\kern-32pt\hbox{$\backslash$}\;\;\;\; v_{k+1}^2 dx\right)^\frac 12\leq \frac 12\left(\int_{B_{\frac 1{2^{k}}}}\kern-22pt\hbox{$\backslash$}\;\; \;\; v_k^2 dx\right)^\frac 12.
$$
Now from (\ref{key1}), we have 
$$m_0+m_1+...\leq m_0+\frac 12 m_0+\frac 1{2^2}m_0+\cdots=2m_0.$$
Noticing that $v_{k}=(u-A(m_0+m_1+...+m_{k-1}))^+\geq (u-2Am_0)^+$,  we have
$$
\left(\displaystyle\int_{B_{\frac 1{2^{k}}}}\kern-22pt\hbox{$\backslash$}\;\; \;\; ((u-2Am_0)^+)^2 dx\right)^\frac 12
\leq 
\left(\displaystyle\int_{B_{\frac 1{2^{k}}}}\kern-22pt\hbox{$\backslash$}\;\; \;\; v_k^2 dx\right)^\frac 12
\leq \frac {m_0}{2^k}\to 0,
$$
which, in conjunction with the assumption that $0$ is a Lebesgue point for $u$,  implies the conclusion.

\vglue 1pc

Now we prove our first step.

\begin{theorem} \label{max}  [Local Maximum Bound]
(a) Let $u$ be a non-negative subsolution of (\ref{Lu}) in $B_2$. 
\noindent Then 
\begin{equation}
\sup\limits_{B_{1}}\ u\le
2^{n+1}A\left(\int_{B_2}\kern-16pt\hbox{$\backslash$} \;\; u^2 dx\right)^\frac 12.\label{inf}
\end{equation}
\end{theorem}
(b) If  $u$ is a nonnegative subsolution in $B_r$, then
\begin{equation}
\sup\limits_{B_{\frac r2}}\ u\le
2^{n+1}A\left(\int_{B_r}\kern-16pt\hbox{$\backslash$} \;\; u^2 dx\right)^\frac 12.
\label{supr}
\end{equation}
(c) If $u$ is a weak solution in $B_r$ then
\begin{equation}
\sup\limits_{B_{\frac r2}}\ |u|\le
2^{n+1}A\left(\int_{B_r}\kern-16pt\hbox{$\backslash$} \;\;|u|^2 dx\right)^\frac 12.\label{maxin}
\end{equation}

\vglue 1pc

{\bf Proof.}
(a) We apply Lemma \ref{point-infty} to $u(x+y)$ for each Lebesgue point $y\in B_1$ of $u$ to obtain
$$
u(y)\leq 2A\left(\int_{B_1}\kern-16pt\hbox{$\backslash$} \;\; (u(x+y))^2 dx\right)^\frac 12\leq 2^{n+1}A\left(\int_{B_2}\kern-16pt\hbox{$\backslash$} \;\; u^2 dx\right)^\frac 12.
$$
Now we are done by taking supremum in all such $y$'s.

(b) Follows from (a) by scaling.

(c)  Notice that both $u^+$ and $u^-$ are nonnegative subsolutions and the conclusion follows from (b) immediately. 

\vglue 1pc

Now we are reached our first milestone of our proof.

\section{Oscillation Estimates}
Our next task is to show the continuity of the solution for which we have to lower the  upper bound or raise the lower bound from $B_2$ to $B_1$.
We start our with a situation which we can lower  the upper bound.
\vglue 1pc

\begin{lemma}\label{smallm}
There is a constant $\epsilon_0>0$ so that if  $u $ is a subsolution in $B_1$ and that $0\leq u\leq 1$ and that 
$$
|\{x\in B_1:u> 0\}|\leq \epsilon_0|B_1|, i.e. |\{x\in B_1:u(x)=0\}|\geq (1-\epsilon_0)|B_1|
$$
then  for $x\in B_\frac 12$, 
$$
u(x)\leq \frac 12.
$$
\end{lemma}
{\bf Proof.}
We apply to the local maximum bound (\ref{supr}) to the subsolution $u$
$$
\sup_{B_\frac 12} u\leq \frac{2^{n+1}A}{|B_1|^\frac 12}\left(\int_{B_1} u^2\right)^\frac 12\leq \frac{2^{n+1}A}{|B_1|^\frac 12}\left|\{ u> 0\}\cap B_1\right|^\frac 12\leq
{2^{n+1}A}{\epsilon_0}^\frac 12= \frac 12,
$$
if $\epsilon_0=\frac{1}{{2^{2n+4}A^2}}$ and that is $u\leq \frac 12$ in $B_\frac 12$.

\vglue 1pc

Our next two lemmas are related with the so called the shooting method. The first is shooting from one point to the sphere and the second is about shooting from a set to another, which we will use to show that the solution needs some space to go up from a levelset to another as set forth in Lemma \ref{key2}.

\begin{lemma}\label {l2} 
Suppose $E$ is a measurable subset of $B_1$ and  $x\in B_1$. Let $\Sigma(E, x)=\{\sigma\in S^{n-1}=\partial B_1:x+ r \sigma\in E\ \text{for some }\ r>0\}$, i.e, $\Sigma(E, x)$ are the collection of the directions in which $E$ can be observed from $x$.

Then 
$$|\Sigma(E, x)|\ge \frac{|E|}{2^{n}}.$$

\end{lemma} 

\vglue 1pc

{\bf Proof.} We use polar coordinates with center $x$ as
\begin{align*} 
|E|&=\int_0^2|(x+r S^{n-1})\cap E| dr\\
&\le 2|(x+r_0S^{n-1})\cap E|\\
&\le 2^{n}|\Sigma(E, x)|\ ,
\end{align*}
\noindent by the mean value theorem and  that for every point $y\in (x+r_0S^{n-1})\cap E$, $y=x+r_0\sigma$ and therefore $\sigma\in \Sigma(E, x)$.
This completes the proof.

\vglue 1pc

Our second shooting lemma shows that one can always find a direction so that the shadow of a measurable set always covers a part of another measurable set in a precise fashion. This shall have of independent interests.

\begin{lemma}\label{shooting}
Suppose $E_1, E_2$ are disjoint measurable subsets of $B_1$.
Then there are $\sigma_0\in S^{n-1}$ and $E_3\subset E_2$ with 
\begin{equation}\label{e3}
|E_3|\geq \frac{|E_1||E_2|}{2^{n}| S^{n-1}|}
\end{equation} 
so that for each $x\in E_3$ there is a $t>0$ such that $x+t\sigma_0\in E_1$, i.e. a part of $E_2$ is covered by $E_3$ in the viewing direction $\sigma_0$.
\end{lemma}

{\bf Proof.}
Let 
$$
\Sigma=\{(\sigma, x)\in S^{n-1}\times B_1: \sigma\in \Sigma(E_1, x), x\in E_2\}.
$$
Then we integrate its indicator function and compute its measure as
$$
\int_{S^{n-1}}{\scalebox{1.5} {$\chi$}}_\Sigma(\sigma, x)d\sigma=|\Sigma(E_1, x)|\geq \frac {|E_1|}{2^{n}}.
$$
Integrate in $x\in E_2$ and compute as
$$
\int_{S^{n-1}\times B_1}{\scalebox{1.5} {$\chi$}}_\Sigma(\sigma, x)d\sigma dx=\int_{S^{n-1}\times E_2}\scalebox{1.5} {$\chi$}_\Sigma(\sigma, x)d\sigma dx\geq \frac {|E_1||E_2|}{2^{n}}.
$$
We change the order of integration and apply the mean value theorem to obtain

$$\int_{S^{n-1}\times B_1}\scalebox{1.5} {$\chi$}_\Sigma(\sigma, x)d\sigma dx
\leq |S^{n-1}|\int_{E_2}\scalebox{1.5} {$\chi$}_\Sigma(\sigma_0, x)dx,
$$
for some $\sigma_0\in S^{n-1}$. 
Now we have
$$\int_{E_2}\scalebox{1.5} {$\chi$}_\Sigma(\sigma_0, x)dx=|\{x\in E_2: \sigma_0\in \Sigma(E_1, x)\}|
\geq \frac {|E_1||E_2|}{2^{n}| S^{n-1}|}.$$
The theorem follows by setting
$$
E_3=\{x\in E_2: \sigma_0\in \Sigma(E_1, x)\}.
$$

{\bf Remark 3.} Our proof actual yields a stronger result. We define
$$
S(E_1)=\inf\{|\Sigma(E_1, x)|:x\in B_1\},
$$
which is geometrically the minimum shadow of $E_1$ projected from a point in the ball.
Then the proof above shows that
\begin{equation}\label{shadow}
|E_3|\geq \frac{S(E_1)|E_2|}{| S^{n-1}|}.
\end{equation} 
We shall that the estimate (\ref{shadow}) is  scaling invariant and also optimal shown in the following example. The shadow is a kind of a measure of its projection but it is quite different from the Besicovitch-Federer's structure theorem \cite {Federer}, where one has to go to all the possible subsets but here one sees only the boundary as one can see even if the tree at $E_1$ is hollow.  We show that (\ref{shadow}) is optimal in the following example.

Image we are at the north pole where the sun is shinning at the horizon around it. $E_1$ a big tree located at the $(0, \frac 78)$ on the $y-$axis with diameter $\frac \epsilon 2$ and $E_2$ are the row of trees on the $x-$axis aligned the interval $[-1, -\frac 12]$ while each tree has diameter $\delta$ and the trees are $\epsilon$ distance apart with total about $\frac 1{2\epsilon}$ many trees there. 
Therefore, $S(E_1)\sim\epsilon$ and $|E_2|\sim\frac{\delta^2}\epsilon$. Then $E_3$ is the set when the shadow of the trees in $E_2$ lands on $E_1$. It is easy to see that $|E_3|= \pi\delta^2\approx S(E_1)|E_2|$.

\vspace{0.5in}

\begin{figure}
\centering
\hspace{-2.0in}\includegraphics[scale=0.5]{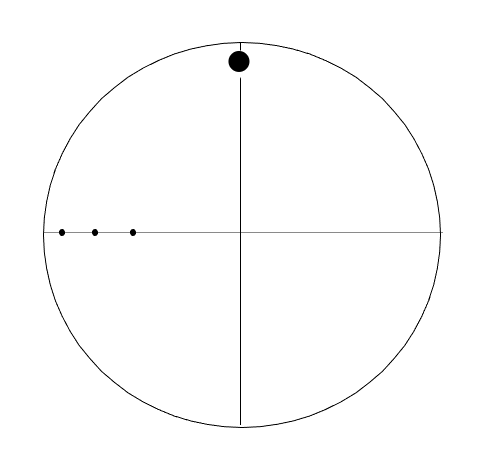}
\end{figure}

\vglue 1pc

Our next step is to find a situation so that Lemma \ref{smallm} can be applied.

\begin{lemma}\label{key2}
Suppose $u$ is a subsolution in $B_2$ with $0\leq u\leq 1$ and let $\kappa>0$
and also suppose
 $|\{x\in B_1: u=0\}|
 \geq \frac {|B_1|}{2}$ and that 
 $|\{x\in B_1: u\geq  \frac 12\}|\geq \kappa|B_1|$,
then 
$$
|\{x\in B_1:0<u(x)<\frac 12\} |\geq \frac {\kappa^2}{n^22^{3n+6}C_1|B_1|}.
$$
\end{lemma}
{\bf Remark 3.}
The geometry of this lemma is the inverse of $u$,
$$
u^{-1}(0, \frac 12):=\{x\in B_1: 0<u(x)<\frac 12\}
$$ must has certain measure if $u$ is of finite energy and that the volume of the levels that these level jumps over are nonzero. It can be viewed as a version of solving an equation and counting the solutions by its measure.

{\bf Proof.} First by the energy estimate (\ref{energyr}) on $u$ with $r=2$, we have that
$$
\int_{B_1}|Du|^2 dx\leq \frac {C_1}4\int_{B_2}u^2\leq  2^{n-2}C_1|B_1|.
$$
We apply Lemma \ref{shooting},  on $E_1=\{x\in B_1:u=0\}$ and $E_2=\{x\in B_1:u\geq \frac 12\}$,
to find  a direction $\sigma_0\in S^{n-1}$ so that 
$$
|E_3|=|\{x\in E_2: \sigma_0\in \Sigma(E_1, x)\}|\geq \frac{|E_1||E_2|}{2^n|\partial B_1|}\geq \frac {\kappa|B_1|}{n2^{n+1}}.
$$
Without loss of generality, we suppose $\sigma_0=e_1=(1, 0, \dots, 0)$.

Now for each $x\in E_3$,   there is  a $T(x)>0$  so that $y=x+T(x)e_1\in E_1$ and thus
 $u(x)\geq \frac 12,  u(y)=0$. With  Fubini theorem, we have that for almost every $x\in E_3$,  we have that $u$ is absolute continuous along  $L_x$ and is in $H^1(L_x)\subset C^\frac 12(L_x)$ where  the line segment $Lx=\{x+te_1: t\in R^1\}\cap B_1$.  Consequently  we can find   $0\leq t_1(x)<t_2(x)\leq T(x)$ so that 
$$
\frac 12=u(x+t_1(x)e_1)>u (x +te_1)>u(x+t_2(x)e_1)=0,
$$ 
for all $t_1(x)<t<t_2(x)$. This is saying that the line segment:
$$
\{x+te_1: t_1(x)<t<t_2(x)\}\subset \{x\in B_1: 0<u(x)<\frac 12\}.
$$
Thus 
$$
\frac 12=u(x+t_1(x)e_1)-u(x+t_2(x)e_1)=\int_{t_2(x)}^{t_1(x)} {D_1}(x+te_1)dt,
$$
Now integrate on  $P_{e_1^\perp}(E_3)=\{(x_2, \cdots, x_n): (x_1, x_2, \cdots, x_n)\in E_3\}$, the project of $E_3$ to the hyperplane $e_1^\perp$, 
\begin{align*}
\frac {|P_{e_1^\perp}(E_3)|}2&\leq \int_{P_{e_1^\perp}(E_3)}\int_{t_2(x)}^{t_1(x)} D_1u(x+te_1)dtdx_2\cdots dx_n\\
&\leq \int_{\{x\in B_1: 0<u<\frac 12\}}|Du|dx\\
&\leq
|\{x\in B_1: 0<u<\frac 12\}|^\frac 12(\int_{B_1}|Du|^2dx)^\frac 12\\
&\leq (2^{-1}C_1|B_1|)^\frac 12 |\{x\in B_1: 0<u<\frac 12\}|^\frac 12.
\end{align*}

Notice that $E_3\subset P_{e_1^\perp}(E_3)\times (-1, 1)$, and therefore
$$
\frac {\kappa}{n2^{n+3}}\leq \frac{|E_3|}4\leq\frac {|P_{e_1^\perp}(E_3)|}2\leq \frac 12(2^{n-2}C_1|B_1|)^\frac 12 |\{x\in B_1: 0<u<\frac 12\}|^\frac 12, $$
and finally 
$$
|\{x\in B_1: 0<u<\frac 12\}|\geq \frac {\kappa^2}{n^22^{3n+2}C_1|B_1|}. 
$$
\vglue 1pc


\begin{theorem}\label{osc}
There is $0<\gamma<1$ depending only on ellipticity constant and dimension $n$ such that if $u$ is a solution in $B_2$,
then
$$
\osc_{B_\frac 12} u\leq \gamma \osc_{B_2} u
$$
where $\osc_{B_r}u=\sup_{B_r} u-\inf_{B_r}u$, the oscillation of $u$ in $B_r$.
\end{theorem}

 {\bf Proof.}  We may assume that $0\leq u\leq 1$ in $B_2$ and further by considering either
 $1-u$ or
$u$  we may also assume that 
$|\{x\in B_1: u\leq \frac 12 \}|\geq \frac 12 |B_1|.$

We want to show $u(x)\le \gamma<1$ for all $x\in B_1$ with some $\gamma$ universal.  

Now examine the measure
$$
|\{x\in B_1: u\geq \frac 34\}|=:\kappa |B_1|.
$$
Now if $\kappa\leq \epsilon_0$, defined in Lemma\ref{smallm}, then we 
$$
u_1(x)=2(u-\frac 12)^+
$$
satisfies the conditions of Lemma \ref{smallm} and thus $u_1\leq \frac 12$ in $B_\frac 12$ therefore we have
$u\leq \frac 34$ and we are done by taking $\gamma=\frac 34$.

Otherwise, it is $\kappa\geq \epsilon_0$, then we can apply Lemma \ref{key2}, to obtain that
\begin{equation}
|\{x\in B_1: 0<u_1(x)<\frac 12\}|\geq \frac {\kappa^2}{{n^2 }{2^{3n+2}}C_1|B_1|}=:\epsilon_1|B_1|.
\end{equation}
This leads to 
\begin{align*}
|\{x\in B_1: u\leq \frac 34\}|&= |\{x\in B_1:u(x)\leq \frac 12\}\cup \{x\in B_1:\frac 12<u(x)\leq \frac 34\}|\\
&\geq |\{x\in B_1:u_1(x)=0\}|+|\{x\in B_1:0<u_1(x)<\frac 12\}|\\
&\geq (\frac 12+\epsilon_1)|B_1|.
\end{align*}
This process can be repeated by considering
$$
u_k=2^k(u-(1-2^{-k}))^+
$$
\noindent and has to stop before $\{x\in B_1: u\geq 1-\frac 1{2^k}\}$ reaches the situation with  measure density less than $\epsilon_0$ in $B_1$ before $k=k_0= \frac {|B_1|}{2\epsilon_1}$ steps and thus 
$$
u\leq 1-\frac 1{2^{k_0+1}}
$$
in $B_\frac 12$.

\begin{theorem}
If $u$ is a weak solution of (\ref{Lu}) in $B_2$, then $u$ is H\"older continuous in $B_1$ with estimates.
\end{theorem}
{\bf Proof:} 
This follows from Theorem \ref{max} and \ref{osc} using the standard iterations.

\vglue 1pc

{\bf Acknowledgement} This research is support in part by a grant from Simons Foundation.

\begin{tabular}{ll}
Department of Mathematics&School of Mathematical Sciences \\
University of Iowa&Shanghai Jiaotong University\\
Iowa City, IA52242&Shanghai, 200240\end{tabular}
\end{document}